%% file: main.tex
\documentclass{amsart}
\input{header}

\DeclareMathOperator{\drop}{drop}

\newcommand{\mcolor}{red}
\newcommand{\msays}[1]{{\footnotesize\textcolor{\mcolor}{\textbf{M:} #1}}}
\newcommand{\tcolor}{blue}
\newcommand{\tsays}[1]{{\footnotesize\textcolor{\tcolor}{\textbf{T:} #1}}}

\newcommand{\TODO}{\textcolor{red}{\footnotesize TODO}}

\begin{document}

	
\title{Rational barycentric coordinates and low degree adjoints}
		
\author{T.\ Baumbach\,}
\author{\,M.\ Winter}
	
%
%
		
\date{\today}
\begin{abstract}
We study rational generalized barycentric coordinates (rational GBCs) on polytopes.
We recall the their essential theory and clarify some statements from the literature.
We study polytopes that admit GBCs with low degree denominators (and hence, low degree (affine) adjoint polynomials).
We give a full classification of polytopes that admit polynomial barycentric coordinates (\ie\ constant adjoints).
\end{abstract}





\maketitle


\input{sec/Introduction}









\par\bigskip
\noindent
\subsection*{Statement on AI use}
AI was used only for finding several of the examples in \cref{sec:appendix_examples}.
The article is entirely human written and AI had no part in the presented research results outside the appendix.

\par\bigskip
\noindent
\subsection*{Funding}
Martin Winter is funded by the SPP 2458 ``Combinatorial Synergies'' (project ID 539851419), funded by the Deutsche Forschungsgemeinschaft (DFG, German Research Foundation).


\bibliographystyle{abbrv}
\bibliography{literature}

\newpage
\input{sec/appendix}

\end{document}

%% file: header.tex
\usepackage[T1]{fontenc}
\usepackage[utf8]{inputenc}
\usepackage[USenglish]{babel}
\usepackage{comment}
\usepackage[
	bookmarks=true,
	plainpages=false,
	linktocpage,
	colorlinks=true,
	citecolor=green!80!black,
	linkcolor=red!70!black,
	filecolor=magenta,
	urlcolor=magenta,
	breaklinks,
	pdfauthor={Martin Winter},
]{hyperref}

\usepackage{amsmath,amsthm}
\usepackage{calc, mathtools} 

\iftrue
\usepackage{amssymb} 
\else
\usepackage[charter,cal=cmcal]{mathdesign}
\fi

\usepackage[
    colorinlistoftodos,
    backgroundcolor=yellow,
    textsize=footnotesize 
    , disable
]{todonotes}
\newcommand{\notodo}[1]{\todo[backgroundcolor=white]{#1}}

\usepackage{stmaryrd} 

\usepackage{colortbl,color} 
\usepackage[dvipsnames]{xcolor}
\usepackage{centernot} 
\usepackage{array} 
\usepackage{enumitem,moreenum} 
\usepackage[noadjust]{cite} 
\usepackage[nameinlink,capitalize,noabbrev]{cleveref}
\usepackage{nicefrac}

\usepackage{tikz-cd} 

\usepackage[font=small,labelfont=bf]{caption}

\usepackage{blkarray}

\usepackage{inconsolata}

\usepackage{accents}
\newlength{\dhatheight}

\definecolor{NiceBlue}{rgb}{0.15, 0.2, 0.75}
\newcommand{\Def}[1]{\textcolor{NiceBlue!90}{\textit{#1}}}

\newcommand{\RR}{\mathbb{R}}    
                   
\newcommand{\bs}[1]{\boldsymbol{#1}}

\newcommand{\tildebs}[1]{\tilde{\boldsymbol{{#1}}}}
\newcommand{\tilbs}[1]{\tildebs{#1}}
\newcommand{\hatbs}[1]{\hat{\boldsymbol{#1}}}

\newcommand{\hatij}{{\ensuremath{\kern1pt\widehat{\kern-1pt\imath\kern-2.3pt\jmath}}}}

\def\:{\colon}

\newcommand{\labelstyle}[1]{\upshape(\textit{#1})}
\newcommand{\mylabel}{\labelstyle{\roman*}}
\newenvironment{myenumerate}{\begin{enumerate}[label=\mylabel]}{\end{enumerate}}

\def\itm#1{{\labelstyle{\romannumeral#1\relax}}}

\def\nlspace{\nolinebreak\space}
\def\nls{\nlspace}

\newcommand{\freespace}{\kern.07em}

\newcommand{\ul}[1]{\underline{\smash{#1}}}

\numberwithin{equation}{section}

\newtheoremstyle{mythmstyle} 
    {\parsep}                    
    {\parsep}                    
    {\itshape}                   
    {}                           
    {\bfseries\scshape}          
    {.}                          
    {.5em}                       
    {}  
\newtheoremstyle{mydefstyle} 
    {\parsep}                    
    {\parsep}                    
    {}                   
    {}                           
    {\mdseries\scshape}          
    {.}                          
    {.5em}                       
    {}  

\theoremstyle{theorem}
\newtheorem{theorem}{Theorem}[section]
\newtheorem{corollary}[theorem]{Corollary}
\newtheorem{lemma}[theorem]{Lemma}
\newtheorem{proposition}[theorem]{Proposition}
\newtheorem{conjecture}[theorem]{Conjecture}

\theoremstyle{definition}

\newtheorem{example}[theorem]{Example}
\newtheorem{remark}[theorem]{Remark}
\newtheorem{question}[theorem]{Question}

\crefname{theorem}{Theorem}{Theorems}
\crefname{proposition}{Proposition}{Propositions}
\crefname{lemma}{Lemma}{Lemmas}
\crefname{corollary}{Corollary}{Corollaries}
\crefname{remark}{Remark}{Remarks}
\crefname{example}{Example}{Examples}
\crefname{definition}{Definition}{Definitions}
\crefname{problem}{Problem}{Problems}
\crefname{observation}{Observation}{Observation}
\crefname{construction}{Construction}{Construction}

\theoremstyle{theorem}
\newtheorem{innercustomgeneric}{\customgenericname}
\providecommand{\customgenericname}{}
\newcommand{\newcustomtheorem}[2]{%
  \newenvironment{#1}[1]
  {%
   \renewcommand\customgenericname{#2}%
   \renewcommand\theinnercustomgeneric{##1}%
   \innercustomgeneric
  }
  {\endinnercustomgeneric}
}

\newcustomtheorem{theoremX}{Theorem}
\newcustomtheorem{lemmaX}{Lemma}
\newcustomtheorem{conjectureX}{Conjecture}

\DeclareMathOperator{\aff}{aff}
\DeclareMathOperator{\conv}{conv}

\DeclareMathOperator{\vol}{vol}  	
\DeclareMathOperator{\Int}{int}

\DeclareMathOperator{\ex}{ex}   
\DeclareMathOperator{\Fac}{Fac}

\DeclareMathOperator{\adj}{adj}

\DeclareMathOperator{\codim}{codim}

\let\eset=\varnothing

\let\<=\langle
\let\>=\rangle
\let\x=\times

\def\...{...}
\newcommand{\shortStyle}{\textit}
\newcommand{\ie}{\shortStyle{i.e.,}}

\newcommand{\eg}{\shortStyle{e.g.}}

\newcommand{\cf}{\shortStyle{cf.}}

\newcommand{\resp}{resp.}

\makeatletter
\renewcommand*{\eqref}[1]{%
  \hyperref[{#1}]{\textup{\tagform@{\ref*{#1}}}}%
}
\makeatother



%% file: sec/Introduction.tex
\newcommand{\p}{p}
\newcommand{\q}{q}

\section{Introduction}

Given a convex polytope $P\subset\RR^d$, every point $x\in P$ can be expressed as a convex combination of the polytopes vertices.
If $P$ is a simplex, then this representation~is unique and the convex coefficients are known as the \emph{barycentric coordinates} of~$x$.\nls
In contrast, if $P$ has more than $d+1$ vertices, there is no longer a unique~such representation, but there are many.
This motivates the following definition:
%
a system of \Def{generalized barycentric coordinates} on $P$ (\Def{GBCs} for short) consists of a family of maps $\bs\alpha=\{\alpha_v:P\to\RR\}$, one per vertex $v$ of $P$, so that for all $x\in P$ holds
%
\vspace{0.5ex}
\begin{align*}
    \text{(\Def{positivity})} \quad\qquad\,\, \alpha_v(x) &\ge 0,
    \\[1.3ex]
    \text{(\Def{normalization})} \quad\;\,\kern0.75pt \sum_v \alpha_v(x) &= 1,
    \\
    \text{\raisebox{0.2ex}{(\Def{precision})}} \quad \sum_v \alpha_v(x) v &= x.
\end{align*}

All sums range over the vertices $v$ of $P$.
In other words, GBCs make~a~consistent choice of convex representation across all points of the polytope.

{\color{black}

GBCs provide a solution to the ubiquitous problem of interpolating data given only at the vertices of a polytope, as it appears, for example, in finite element~analysis, computer graphics or geometric modeling.
An overview of different well-studied GBCs and some applications can be found in \cite{floater2015generalized}.
For many of these applications it is advantageous to use GBCs that are easy to represent and evaluate. This constraint naturally favors algebraic over transcendental expressions and motivates the study of polynomial and rational GBCs.

In this paper we study such polynomial and rational GBCs, with a particular~focus on understanding their degrees.


\subsection{Polynomial GBCs}
\label{sec:polynomial_GBCs}

\Def{Polynomial GBCs} have $\alpha_v\in\RR[x]$ a multivariate polynomial in $x\in\RR^d$.
The classical barycentric coordinates on simplices are linear,\nls and hence examples thereof.
Other examples are obtained by taking products of polytopes and GBCs (\cf\ \cref{ex:product}).



As pointed out by Wachspress already in 1975, most polytopes do not admit polynomial GBCs:

\begin{theorem}[{\cite[Chapter 2, §2.1]{wachspress1975rational}}]
    A quadrangle has polynomial GBCs if and only if it is a parallelogram.
\end{theorem}

The polynomial GBCs on a parallelogram are obtain as affine transformations of the product GBCs on squares.
We likewise obtain polynomial GBCs on all affine transformations of products of simplices.
In \cref{sec:polynomial_GBCs} we prove that this list is complete:
%
%
%

\begin{theorem}
    \label{res:polynomial_GBCs}
    A polytope admits polynomial GBCs if and only if it is affinely~equi\-valent to a product of simplices.
\end{theorem}

This resolves \cite[Conjecture 4.8]{clarke2020moment}.

\subsection{Rational GBCs}

In his work on the finite element method Wachspress gave a construction for \Def{rational GBCs} that exist for all convex polygons and simple 3-polytopes \cite{wachspress1975rational}
%
%
(\ie\ $\alpha_v=\p_v/\q$ where $\p_v,\q\in\RR[x]$ are multivariate polynomials in $x\in\RR^d$).
Warren later extended the construction to polytopes of general dimension and combinatorics \cite{warren1996barycentric}.
These coordinates came to be known as the \emph{Wachspress coordinates} and are of great practical and theoretical interest (\cf\ \cref{ex:Wachspress}).

We emphasize that the typical definition of rational GBCs (and a main reason for the nice theory) is that all coordinates $p_v/q$ are over \emph{a common} denominator $q$.
We will assume throughout that $q$ does not vanish on the interior of the polytope.






In \cref{sec:Warren} and \cref{sec:min_numerator_degree} we recall and complement Warren's theory of denominator degrees in rational GBCs.
We provide corrections and completions to important points in Warren's arguments.

In \cref{sec:low_denominator_degrees} we transition to studying the denominator degrees of rational GBCs. We in particular focus on low degree denominators.
In particular, the denominator degree is zero if and only if the GBCs are polynomial.
In \cref{sec:polynomial_GBCs} we apply this theory to provide the classification of polytopes with polynomial GBCs (\cref{res:polynomial_GBCs}).
This amounts to the classification of polytopes with degree-zero adjoints.

\section{On Warren's theory of rational GBCs}
\label{sec:Warren}
\label{sec:structure}

This section focuses primarily on the theory of rational GBCs developed by Warren, and in particular, on reproving degree bounds and the uniqueness of Wachspress coordinates.
For this we essentially follow \cite{warren2003uniqueness}.
At times this section goes beyond the classical work and provides deeper insight into the structure of degrees and vanishing orders of general rational GBCs.
%

\begin{remark}

    While \cite{warren2003uniqueness} contains the essential ideas, it also contains a number of inaccuracies and gaps that have confused readers over the years.
    Some of them are not trivial to patch.
    A few stem from being imprecise about whether the setting~is affine or projective, which leads to neglecting the consequences of parallelisms in~the face structure of polytopes.
    In fact, there are true subtleties to the degree theory once we work in affine space (\eg\ the denominator degree is not bounded from below by $\#\text{facets}-\text{dim}-1$).
    %
    Other gaps appear in central proofs, such as local-global arguments for proving degree bounds.
    We include comments throughout the section (in footnotes and remarks) to point out whenever we repair or elaborate on \cite{warren2003uniqueness}.\todo{\msays{Are the footnotes there?}}


    
\end{remark}

Throughout, let $P\subset\RR^d$ be a $d$-dimensional polytope (\ie\ with non-empty~interior) and $\alpha_v=\{p_v/q\}$ rational GBCs on $P$.


\subsection{Affine precision}

On rational GBCs \textit{normalization} and \textit{precision} take on the following form:
\begin{equation}
    \label{eq:q_sum_pv}
    q = \sum_v p_v \quad\text{and}\quad  xq = \sum_v p_v v.
\end{equation}
The left sum states that $q$ is determined by the numerator polynomials.
%
A standard argument shows that the pair \eqref{eq:q_sum_pv} is equivalent to the single property called \Def{affine precision}: for any affine linear form $L$ holds
\begin{equation}
    \label{ex:affine_precision}
    L(x) q(x) = \sum_v L(v) p_v(x).
\end{equation}
Some authors define GBCs (rational or not) in terms of this single property. 

\subsection{Positivity}
\label{sec:positivity}

The third defining property, \emph{positivity}, is used essentially in the~following way:
%
%
if $F\subset P$ is a facet and $v\in P$ is a vertex that does not lie in $F$, then we necessarily have $\alpha_v(x)=0$ for all $x\in F$, or equivalently, $p_v(x)=0$ for all $x\in F$. 
In still other words, the \Def{facet linear form $\ell_F(x)$} (characterized by~$F = P \cap \{\ell_F(x)=0\}$) divides the numerator polynomial $p_v$.
We conclude that the numerators of general rational GBCs factors as follows:
\begin{equation}
    \label{eq:factorization}
    p_v(x)=\hat p_v(x) \prod_{F\not\ni v} \ell_F(x),
\end{equation}
where the product is over all facets $F\subset P$ that do not contain the vertex $v$, and~$\hat p_v$ is some polynomial whose form we do not yet know.

Interestingly, \eqref{eq:factorization} is the only consequence of positivity that we will make use of.
Hence, at least for our purpose, one could replace positivity by \eqref{eq:factorization}. 


\subsection{Degrees}
\label{sec:degrees}

There are essentially three ways to measure complexity of GBCs by their degrees: numerator degrees, denominator degrees and rational degrees.


\emph{Numerator degrees} have been the primary focus of Warren's work.
For example, consider GBCs on a simple polytope $P$.
Then from \eqref{eq:factorization} follows that the degree of $p_v$ must be at least $m-d$, where $m$ is the number of facets of $P$.
There are three major insights from Warren's theory: 
\begin{enumerate}[label=(\arabic*)]
    \item the lower bound of $m-d$ on the numerator degree applies to all polytopes, not just simple polytopes (we reprove this in \cref{res:lower_bound_degree_GBCs});
    \item this lower bound is actually attained by some GBCs
    (see \cref{ex:Wachspress});
    \item the rational GBCs that attain this lower bound are essentially unique (we prove this in \cref{res:GBC_deg_m_d_is_WPC}).
\end{enumerate}
These unique degree minimizers are known as the \Def{Wachspress coordinates $\bs\alpha^{\mathrm W}$} of~$P$ and have many interesting properties and use cases.
We emphasize again that this extraordinarily nice theory is, to a decent extent, a consequence of the fact that all coordinates of GBCs share the same denominator (see also \cref{ex:pyramid} and \cref{ex:joins}).


\emph{Denominator degrees}, in contrast, are much less  well-understood.
The right sum in \eqref{eq:q_sum_pv} gives
\begin{equation}
    \label{eq:deg_q_le_deg_p_1}
    \deg(q) \le \max_v \deg(p_v) - 1
\end{equation}

The denominator of Wachspress coordinates is called the \Def{adjoint} or \Def{adjoint polynomial $\adj_P$} of $P$.
There are several ways to define this polynomial even without~the Wachspress coordinates.
This often happens in a projective setting, where one can show that the (homogeneous) degree is exactly $\deg(\adj_P)=m-d-1$ \cite{warren1996barycentric}.
One of the oversights of Warren in \cite{uniqueCentralSymmetry} was that this does not carry over to the affine setting.
For a general polytope, the degree of $\adj_P$ can be much smaller than $m-d-1$.
And so, while for a generic polytope the Wachspress coordinates satisfy \eqref{eq:deg_q_le_deg_p_1} with equality, this is not the case for general polytopes, not even if they are simple or simplicial.
This ``deficiency in degree'' 
$$\drop(P):=m-d-1-\deg(\adj_P)$$ 
was studied in detail in \cite{baumbach2025canonicalformscissorscongruence}, where it was found to subtly relate to parallelisms~in~the face structure of $P$.

For general GBCs, no nice theory (as compared with the numerator degrees) has emerged so far. In contrast, we found polytopes $P$ with GBCs $\{p_v/q\}$ so that
\begin{enumerate}[label=(\arabic*)]
    \item $\deg(q)=\deg(\adj_P)$, where $\{p_v/q\}$ are not the Wachspress coordinates.
    \item $\deg(q)<\deg(\adj_P)$.
\end{enumerate}
In \cref{sec:low_denominator_degrees} we introduce $\deg(P)$ as the minimally possible denominator degree of rational GBCs on $P$.

Lastly, one can study the \emph{rational degree} $\deg(\alpha_v):=\deg(p_v)-\deg(q)$.
By \eqref{eq:deg_q_le_deg_p_1} $\deg(\alpha_v)\ge 1$.
We did not pursue this route.
We note that in the case of Wachspress coordinates, $\deg(\alpha^{\mathrm W}_v) = m-d-\deg(\adj_P)$ is precisely $\drop(P)+1$ and its study is therefore subsumed by \cite{baumbach2025canonicalformscissorscongruence}.


\subsection{Basic examples}
\label{sec:examples}


The classical barycentric coordinates are defined on simplices. They are linear functions and are the unique GBCs with this property.
We give a definition based on cone volumes:

\begin{example}[barycentric coordinates]
    \label{ex:barycentric}
    Let $\Delta$ be a $d$-dimensional simplex.
    Let $F_v\subset \Delta$ be the facet opposite to vertex $v$.
    For a point $x\in\Delta$ let $C_v:=F_v\vee x$ be the cone with base $F_v$ and apex $x$.
    Define $\alpha_v:=\vol(C_v)/\vol(\Delta)$.
    Since the cones $C_v$ partition $\Delta$, we have $\sum_v\alpha_v=1$.
    We also have $\alpha_v(w)=1$ if $v=w$ and $0$ otherwise. Lastly, it is easy to see that $\alpha_v(x)$ is linear in $x$.
    It follows that these are the classical barycentric coordinates on $\Delta$.
\end{example}

Interestingly, the Wachspress coordinates can be defined in terms of cone volumes as well, but they require the polar dual.

\begin{example}[Wachspress coordinates]
    \label{ex:Wachspress}
    Given a polytope $P\subset\RR^d$ with $x\in\Int(P)$.
    Consider its \emph{polar dual} at $x$:
    $$(P-x)^\circ:=\{y\in\RR^d\mid \<y,z-x\>\le 1\text{ for all $z\in P$}\}.$$
    Let $F_v^\diamond$ be the facet of $(P-x)^\circ$ that is dual to the vertex $v\in P$.
    Let $C_v$ be the cone with base $F_v^\diamond$ and apex at the origin.
    Define $\alpha_v:=\vol(C_v)/\vol(P-x)^\circ$.
    These are the Wachspress coordinates of $P$ at $x$.
    The fact that these are indeed GBCs is a consequence of Minkowski's balancing condition.
    This geometric cone volume construction goes back to Ju, Schaefer, Warren and Desbrun \cite{ju2005geometric}.
    Evidently, they are rational.
    We do here not reverify the degree but point to the literature.
\end{example}

An important elementary operations to construct new GBCs from old ones is the product construction:

\begin{example}[product coordinates]
    \label{ex:product}
    Given polytope $P_i\subset \RR^{d_i}$ with $i\in\{1,2\}$ with GBCs $\bs\alpha^i$, one can define the \Def{product coordinates} on $P_1\times P_2$.
    The vertices of $P_1\times P_2$ are of the form $(v_1,v_2)$, where $v_i\in P_i$.
    Then
    $$(\alpha^1\times\alpha^2)_{v_1,v_2} := \alpha^1_{v_1}\cdot \alpha^2_{v_2}.$$
    One can verify that these do indeed give GBCs.
\end{example}

Taking the product of barycentric coordinates on simplices gives polynomial GBCs on products of simplices.
As we will show in \cref{sec:polynomial_GBCs}, these are the only polytopes that admit polynomial GBCs.
The simplest non-simplicial example is the square:

\begin{example}[square]
    \label{ex:square}
    $$P = [0,1]^2 = \conv\{v_{00},v_{10},v_{11},v_{01}\}.$$
    Let $\ell_{ij}$ be the linear form of the facet spanned by $v_i$ and $v_j$. Then
    $$\ell_{00,10}=x_2,\quad \ell_{10,11}=1-x_1,\quad \ell_{11,01}=1-x_2,\quad \ell_{01,00}=x_1.$$
    and
    \begin{align*}
    p_{v_{00}} &= \ell_{10,11}\ell_{11,01},\quad
    p_{v_{10}} = \ell_{11,01}\ell_{01,00},\quad
    \\
    p_{v_{11}} &= \ell_{01,00}\ell_{00,10},\quad
    p_{v_{01}} = \ell_{00,10}\ell_{10,11}.
    \end{align*}
    and
    $q=p_{00}+p_{10}+p_{11}+p_{01} = 2$.
\end{example}


\subsection{Vanishing orders and the face arrangement}

A central tool for proving degree lower bounds for the numerators and denominators is to look at the multiplicities of their zeroes.
Let $o_E(p)$ denote the vanishing order of the polynomial $p$ along a non-empty affine subspace $E$.
Clearly, $\deg(p)\ge o_E(p)$.
For a face $f\subseteq P$~we abbreviate $o_f(p):=o_{\aff(f)}(p)$.

\begin{lemma}
\label{res:order}
Given a $\delta$-face $f\subseteq P$ and vertex $v\in P$:
    \begin{myenumerate}
        \item if $v\in f$ then $o_f(p_v) = o_f(\hat p_v)$. 
        \item $o_f(p_v)\ge o_f(q)$ with equality if and only if $v\in f$.
        %
    \end{myenumerate}
    %
\end{lemma}
\begin{proof}
    Part \itm1: by the factorization \eqref{eq:factorization} 
    we have
    $$o_v(p_v)=o_v(\hat p_v)+o_v\Big(\prod_{F\not\ni v} \ell_F\Big).$$
    But the latter term vanishes since none of the linear forms $\ell_F,F\not\ni v$ is zero on $v$.

    Part \itm2: the limits of $p_v(x)/q(x)$ as $x$ approaches the face $f$ from within $P$ exist and yield GBCs on $f$.
    The limit towards a point in the relative interior of $f$ is zero if and only if $v\not\in f$. \qedhere

\end{proof}

The \Def{face arrangement $\mathfrak F(P)$} of a polytope $P$ is the arrangement of the affine~subspaces generated by face spans (including of $P$ and the empty face) closed under~intersection.
The elements of the arrangement are called \Def{flats}.
The ones generated~by faces are called \Def{facial flats}.
The others are called \Def{residual flats}.
Alternatively, a~flat is facial if and only if it is spanned by the vertices $v\in P$ that is contains.

Face arrangements are only implicit in Warren's treatment of the subject.
They are more visible in his work on the adjoint polynomial \cite{warren1996barycentric}.
The modern treatment~of rational GBCs seems to favor this language, especially in the form of the \emph{residual arrangement} \cite{kohn2020projective}.

\notodo{\msays{Add picture of face arrangement.}}

\begin{lemma}
    \label{res:q_ge_min_pv}
If $E\in\mathfrak F(P)$ is a flat then
    $$o_E(q)\ge \min_{v\not\in E} o_E(p_v) - \begin{cases} 1 & \text{if $E$ is a facial flat} \\ 0 & \text{if $E$ is a residual flat}\end{cases}.$$
\end{lemma}
\begin{proof}
    Let $E'\subset E$ be the affine span of the vertices $v\in P$ that lie in $E$.
    Let $L$~be~a generic affine linear form that vanishes on $E'$.
    Then
    $$o_E(L) + o_E(q) = o_E(Lq)
    \overset{(*)}=o_E\Big(\sum_v L(v) p_v\Big) 
    = o_E\Big(\sum_{v\not\in E} L(v) p_v\Big) \overset{(**)}\ge \min_{v\not\in E} o_E(p_v).$$
    where in $(*)$ we used affine precision, and $(**)$ is the well-known inequality for~the vanishing order of sums of polynomials.
    Since $L$ is linear, $o_E(L)\le 1$.
    If~$E$~is~residual, then it is not spanned by the contained vertices, and so $E'\subsetneq E$. 
    Since $L$ was chosen generic, it does not vanish on $E$, and so $o_E(L)=0$. 
\end{proof}

\subsection{Facial flats}
\label{sec:facial_flats}

Given a flat $E\in\mathfrak F(P)$, its \Def{excess $\ex(E)$} counts how many~more facet hyperplanes of $P$ contain $E$ than one would generically except from only knowing its dimension.
More precisely, the number of facet hyperplanes that contain $E$ is $\codim(E)+\ex(E)$, where $\codim(E):=d-\dim(E)$.

\begin{lemma}
    \label{res:o_q_ge_ex_f}
    Given rational GBCs $\{p_v/q\}$ on $P$ and a face $f\subseteq P$, it holds 
    $$o_f(q)\ge \ex(f).$$
\end{lemma}
\begin{proof}
    We use induction on $\codim(f)$.
    If $\codim(f)=0$, then $f=P$ and $\ex(f)=0$, and the statement is immediate.
    Assume now $\codim(f)\ge 1$.
    
    By \cref{res:q_ge_min_pv} we have $o_f(q)\ge \min_{v\not\in f} o_f(p_v)-1$ and so it suffices to show that $o_f(p_v)\ge \ex(f)+1$ for all $v\not\in f$. Fix such a $v\not\in f$.
    Let $f_v:=f \vee v$~be~the~minimal face~of $P$ that contains both $f$ and $v$. 
    To bound $o_f(p_v)$ we consider two~independent contributions to its value: first, the vanishing order of $p_v$ on the larger face $f_v\supset f$; second, additional contributions from the facet linear forms $\ell_F$ where $F$ contains $f$ but not $v$ (\cf\ the factorization \eqref{eq:factorization}):
    \begin{align}\label{eq:of_pv_ge}
        o_f(p_v) 
        &\ge o_{f_v}(p_v) + \#\{\text{facets of $P$ that contain $f$ but not $v$}\}.
    \end{align}
    The left term we bound using \cref{res:order} \itm2 and the induction hypothesis.
    The latter applies due to $\codim(f_v)<\codim(f)$:
    $$o_{f_v}(p_v) \overset{\mathclap{\text{\itm2}}}= o_{f_v}(q) \overset{\smash{\mathclap{\mathrm{IH}}}}\ge \ex(f_v).$$
    The second term we bound as follows:
    \begin{align*}
        \#\{\text{facets that contain $f$ but not $v$}\} \hspace{-5em}&
        \\ &= \#\{\text{facets that contain $f$ but not $f_v$}\}
        \\&= (\codim(f) + \ex(f)) - (\codim(f_v) +\ex(f_v))
        \\&= (\ex(f) - \ex(f_v)) + (\codim(f) - \codim(f_v)).
        \\&\ge (\ex(f) - \ex(f_v)) + 1,
    \end{align*}
    where in the first equality we use that $f_v$ was chosen minimal; in the second equality we use that the facets that contain $f_v$ are a subset of the facets that contain $f$, that $f$ is contained in $\codim(f)+\ex(f)$ facets, and $f_v$ is contained in $\codim(f_v)+\ex(f_v)$ facets.

    Combining both estimations in \eqref{eq:of_pv_ge} finishes the proof.
    \qedhere

\end{proof}

We obtain our first major insight for numerator degrees:

\begin{theorem}
\label{res:pv_ge_m_d}
\label{res:lower_bound_degree_GBCs}
\label{thm:lower_bound_degree_GBCs}
    Let $P\subset \RR^d$ be a polytope with $m$ facets and $\{p_v/q\}$ rational GBCs on $P$. 
    Then $\deg(p_v) \geq m - d$ for all vertices $v$.
\end{theorem}
\begin{proof}
    Suppose that $v$ lies on $d+k$ facets, that is, $\ex(v)=k$. 
    %
    By \cref{res:order}~\itm1~and \itm2, as well as \cref{res:o_q_ge_ex_f}, we have
    %
    %
    $$\mathllap{(*)\quad }\deg(\hat p_v)\ge o_v(\hat p_v) \overset{\mathclap{\text{\itm1}}}=o_v(p_v)\overset{\mathclap{\text{\itm2}}}=o_v(q)\overset{\smash{\mathclap{\text{\ref{res:o_q_ge_ex_f}}}}}\ge \ex(v)= k.$$
    And hence
    \[
        \deg(p_v) 
        \overset{\eqref{eq:factorization}}= \deg(\hat p_v) + \deg\!\Big(\prod_{F \not\ni v}\ell_F\Big)
        \overset{\smash{(*)}}\ge k + (m - (d + k)) = m - d.\qedhere
    \]
\end{proof}

The lower bound of \cref{res:pv_ge_m_d} should be understood as a consequence of~writing our rational GBCs $p_v/q$ over a common denominator polynomial $q$.
Even though the Wachspress coordinates satisfy $\gcd\{p_v,\adj_P\}=1$, individual $p_v$ may still share factors with the adjoint, and canceling these factors would result in lower numerator degrees.

\begin{example}[Pyramid]
    \label{ex:pyramid}
    Given a polytope $P$ with GBCs $\bs\alpha=\{p_v/q\}$.
    Let $P_\star$~be the pyramid over $P$ with apex $\star$ at the origin. The following gives GBCs $\bs\beta$ of $P_\star$:
    \begin{align*}
        \beta_v(x,t) &= t \alpha_v\Big(\frac x t\Big),\quad\text{for vertices $v\in P$},
        \\
        \beta_{\star}(x,t) &= 1-t.
    \end{align*}
    One can verify that these are indeed rational GBCs of $P_\star$.
    %
    %
    Note that $\deg(\beta_\star)=1$, independent of the combinatorics of $P_\star$. Hence, if $m-d>1$, then by \cref{res:pv_ge_m_d} factors in numerator and denominator must have canceled.

\end{example}

In fact, it is possible that $p_v$ and $\adj_P$ share factors at \emph{all} vertices $v$. 

\begin{example}[Joins]
    \label{ex:joins}
    Let $P$ and $Q$ be polytopes with rational GBCs $\bs\alpha^P$ and~$\bs\alpha^Q$.
    Each point of the \emph{join} $P*Q$ can be written as
    $(1-t)x+ty$ for some $x\in P$, $y\in Q$ and $t\in[0,1]$.
    The following defines GBCs $\bs \beta$ on $P*Q$:
    \begin{align*}
        \beta_v(x,y,t) &= (1-t)\alpha_v^P\Big(\frac{x}{1-t}\Big),\quad\text{for all $v\in P$}
        \\
        \beta_w(x,y,t) &= t\alpha_w^Q\Big(\frac{y}{t}\Big),\quad\text{for all $w\in Q$}
    \end{align*}
    Again, one can verify that these are indeed rational GBCs of $P\star Q$.
    Note that,\nls as written here, the $\beta_v$ only depend on $x$, whereas the $\beta_w$ depend only on $y$.
    Hence, unless $\bs\alpha^P$ and $\bs\beta^P$ had constant denominator, writing $\bs\beta$ over a common denominator requires adding shared factors at all vertices.

\end{example}





The following fact, that GBCs cannot vanish on facets, is of use later.

\begin{lemma}
    \label{res:not_vanish_on_facet}
    Let $F\subset P$ be a facet, and assume either of the following:
    \begin{myenumerate}
        \item $\gcd\{q, p_v\}=1$ for some vertex $v\in P$, or
        \item $\gcd\{p_v\mid v\in P\}=1$,
    \end{myenumerate}
    Then $o_F(q)= o_F(p_w) = 0$ for all vertices $w\in F$.
\end{lemma}
\begin{proof}
    Part \itm1: If $o_F(q)>0$, then $o_F(p_v)> 0$ by \cref{res:order} \itm2. But this means
    $\ell_F\mid q$, as well as $\ell_F \mid p_v$, contradicting the coprimality assumption.
    Consequently, $o_F(q)=0$.
    Via \cref{res:order} \itm2 also follows $o_F(p_w)=0$ for all $w\in F$.

    Part \itm2: As in the proof of \itm1, if $o_F(q)>0$, then $\ell_F\mid p_v$ for all $v$.
    Therefore, $\ell_F\mid \gcd\{p_v\}$.
\end{proof}

Results stating non-vanishing on faces of $P$ will be important later, are however hard to come by in this generality.
\Cref{res:not_vanish_on_facet} is among the strongest claims we can make under current assumptions. 
Even assuming $\gcd\{q,p_v\mid v\in V\}=1$,\nls it~is still possible that $q$ vanishes on a face of codimension two, or that $p_v$ vanishes on a facet $F\not\ni v$ to an order higher than one.

\todo{\msays{Add appendix Example where $q$ (even when reduced) vanishes on ridge (aka vertex of a square). And maybe the $p_v$ counterexample as well. Muss weg}}

\subsection{Residual flats}
\label{sec:residual_flats}

\begin{lemma}
    \label{res:q_vanishes_residual}
    If $P$ is a polytope with rational GBCs $\{p_v/q\}$ and a residual flat $E$, then $o_E(q)\ge 1$.
    In other words, $q$ vanishes on the residual flats of $P$.
\end{lemma}
In the literature this is often phrased as ``$q$ vanishes on the residual arrangement of $P$''.
%
%
\begin{proof}
    \cref{res:q_ge_min_pv} states $o_E(q)\ge \min_{v\not\in E} o_E(p_v)$, thus it suffices to show $o_E(p_v)\ge 1$ for all $v\not\in E$.
    Fix such a $v\not\in E$.
    Since $E$ is cut out by facet hyperplanes, there is some facet $F\subset  P$ whose span contains $E$ but not $v$. Hence $\ell_F\mid p_v$ by \eqref{eq:factorization}, and $o_E(p_v)\ge o_E(\ell_F)=1$.
\end{proof}


%
The following conjecture is a natural common generalization of \cref{res:q_vanishes_residual} and \cref{res:o_q_ge_ex_f}.

\begin{conjecture}
    \label{conj:o_E_q}
    \todo{\msays{I know that this is true (ChatGPT), but the proof is to complicated for us here. Still conjecture?}}%
    If $E\in\mathfrak F(P)$ is a flat, then 
    $$o_E(q)\ge \operatorname{ex}(E) + \begin{cases} 0 & \text{if $E$ is a facial flat} \\ 1 & \text{if $E$ is a residual flat}\end{cases}.$$
\end{conjecture}

In fact, in the case of Wachspress coordinates, where the denominator is the~adjoint polynomial $\adj_P$, these vanishing orders are exact \cite{bruser2025geometry}:
$$o_E(\adj_P) = \operatorname{ex}(E) + \begin{cases} 0 & \text{if $E$ is a facial flat} \\ 1 & \text{if $E$ is a residual flat}\end{cases}.$$


One could expect to prove \cref{conj:o_E_q} by combining ideas from \cref{res:o_q_ge_ex_f} and \cref{res:q_vanishes_residual}.
The proof of the latter relies on the fact that we can estimate~the vanishing orders of the $p_v$.
But so far we know how to do this only if $v\in E$.
In~fact, already an analogue of \cref{res:q_vanishes_residual} for numerators polynomials fails for sufficiently complicated examples 
(see \cref{exp:res_arrangement_non_vanishing}).
Analogous can be proven in special cases, for example, see \cref{res:residual_arrangement}.
\todo{\msays{Add example in appendix. There is a commented out result below that explains why a counterexample must have dimension $\ge 5$. Doktorarbeit}}

\section{Minimal numerator degree}
\label{sec:min_numerator_degree}

In this section we focus specifically on rational GBCs $\{p_v/q\}$ of minimal possible numerator degree, that is, $\deg(p_v)= m-d$ according to \cref{res:lower_bound_degree_GBCs}.
The goal is to show that $\{p_v/q\}$ must be the Wachspress coordinates.
%
%
Again, we broadly follow Warren's work \cite{warren2003uniqueness}, though we make corrections at several points in the proof.

It turns out that quite a number of our statements have a ``per-vertex'' version, that is, need to assume $\deg(p_v)=m-d$ only for a single vertex $v$ rather than~all vertices at the same time.
We make an effort to provide per-vertex versions whenever this is possible and reasonable given our setup.
The following is a first example:



\begin{corollary}
    \label{res:min_coprime}
    \label{res:min_not_vanish_of_facet}
    Suppose $\deg(p_v)=m-d$ for some vertex $v\in P$.
    Then
    \begin{myenumerate}
        \item $\gcd\{p_v\mid v\in P\}=1.$
        \item $o_F(q)=0$ for each facet $F\subset P$. 
        %
    \end{myenumerate}    
\end{corollary}
\begin{proof}
    To prove \itm1 suppose that $r\mid p_w$ for all vertices $w$. 
    Then $r\mid \sum_w p_w=q$ and $r$ can be canceled from all $p_w/q$, yielding GBCs of lower numerator degrees. This contradicts \cref{res:pv_ge_m_d} at $v$.

    To prove \itm2 observe that if $o_F(q)>0$, \cref{res:order} \itm2 implies $o_F(p_w)>0$~for all $w\in F$. 
    From \eqref{eq:factorization} we also get $o_F(p_w)>0$ for all $w\not\in F$.
    In other words, $\ell_F\mid p_w$ for all $w$, contradicting \itm1.
\end{proof}

\subsection{Numerators on residual flats}
\label{sec:numerator_residual}

Recall $\hat p_v$ as define in \eqref{eq:factorization}.



\begin{lemma}
    \label{res:deg_ov_ex}
    \label{res:hatpv_homogeneous}
    If $\deg(p_v)=m-d$ for some vertex $v\in P$, then
    \begin{equation} \label{eq:deg_ov_ex}
        \deg(\hat p_v) 
            = o_v(\hat p_v) = o_v(p_v) 
            = o_v(q)
            = \ex(v).
    \end{equation}
    In particular, $\hat p_v$ is homogeneous at $v$ of degree $\ex(v)$ and is (up to a factor) the~homogeneous part of $p_v$ of lowest degree.\footnote{Consider $\hat p_v$ as a sum of homogeneous parts. Then $\deg(\hat p_v)=\ex(v)$ states that the highest degree term has degree $\ex(v)$, whereas $o_v(\hat p_v)=\ex(v)$ states that the lowest degree term has degree $\ex(v)$. 
    Hence the term of degree $\ex(v)$ is the only term.} 
\end{lemma}
This lemma is later strengthened to general faces (see \cref{res:min_deg_o_f_q_eq_ex}), though that requires more preparation.
\begin{proof}[Proof of \cref{res:deg_ov_ex}]
    From minimal degree and the factorization \eqref{eq:factorization} we get
    \begin{align*}
        m - d \;&= \;\deg(p_v) 
        \\&\overset{\mathclap{\smash{\eqref{eq:factorization}}}}= \;\deg(\hat p_v) + \#\{\text{facets that do not contain $v$}\}
        \\&= \;\deg(\hat p_v) + m - \#\{\text{facets that do contain $v$}\}
        \\&= \;\deg(\hat p_v) + m - (d + \ex(v)),
    \end{align*}
    which rearranges to $(*)\,\deg(\hat p_v)=\ex(v)$.
    Using also \cref{res:order} \itm1 and \itm2, as well as \cref{res:o_q_ge_ex_f} we get the following cyclic  chain of inequalities:
    $$\ex(v)
    \overset{\smash{\mathclap{(*)}}}=\deg(\hat p_v) 
    \ge o_v(\hat p_v) 
    \overset{\smash{\mathclap{\text{\itm1}}}}= o_v(p_v) 
    \overset{\smash{\mathclap{\text{\itm2}}}}= o_v(q)
    \overset{\smash{\mathclap{\ref{res:o_q_ge_ex_f}}}}\ge \ex(v).$$
    Thus all terms must be equal, giving the claimed identities in \eqref{eq:deg_ov_ex}.
\end{proof}

The following is an anlogue of \cref{res:q_vanishes_residual} for numerator polynomials:

\begin{lemma}
    \label{res:residual_arrangement}
    Let $E\in\mathfrak F(P)$ be a residual flat.
    If $\deg(p_v)=m-d$ for some vertex $v\in P$, then
    \begin{myenumerate}
        \item $o_E(p_v)\ge 1$. 
        \item if $v\in E$, then $o_E(\hat p_v)\ge 1$. 
        \notodo{Is this ever used, and do we want to keep it? (it is proved as an intermediate step below)}
    \end{myenumerate}
\end{lemma}
\begin{proof}

    \emph{Case $v\not\in E$}:
    The flat $E$ is the intersection of some facet hyperplanes.
    Thus, $v\not\in F$ for at least one such facets $F$. We get $\ell_{F}\mid p_v$, and $o_E(p_v)\ge o_E(\ell_{F})=1$. 

    \emph{Case $v\in E$}:
    Normalization gives us
    \begin{equation}
        \label{eq:q_eq_pv_plus_sum}
        q = p_v + \sum_{w\ne v} p_w 
    \end{equation}
    %
    %
    This identity must hold in each homogeneous degree (defined at $v$), and so we may restrict to degree $\ex(v)$.
    Write $p^{\ex}$ for the homogeneous part of degree $\ex(v)$~of~a~polynomial $p$.
    By \cref{res:deg_ov_ex} we have $p_v^{\ex}=\hat p_v$.
    Using \cref{res:order} \itm2, for $w\ne v$ we have $o_v(p_w)>o_v(q)=\ex(v)$, and hence $p_w^{\ex}=0$.
    In degree $\ex(v)$ the identity
    \eqref{eq:q_eq_pv_plus_sum} therefore reduces to $q^{\ex}=\hat p_v$.
    Now restrict both sides to $E$.
    Recall that restriction to homogeneous parts and restriction to affine subspaces (through $v$) commutes:
    $$\hat p_v|_E = q^{\ex}|_E = (q|_E)^{\ex} \equiv 0$$
    where $q|_E\equiv 0$ by \cref{res:q_vanishes_residual}.
    We conclude $o_E(p_v) = o_E(\hat p_v) \ge 1$.
\end{proof}


\begin{remark}
     \cref{res:deg_ov_ex} combined with \cref{res:residual_arrangement} \itm2 states that $\hat p_v$ is a homo\-geneous polynomial of minimal degree that vanishes on the residual flats of the~tangent cone at $v$.
     The reader familiar with the literature on adjoints might then~recog\-nize that $\hat p_v$ must be the adjoint of the tangent cone \cite{kohn2020projective,bruser2025geometry}.
     This insight could shortcut our proof of uniqueness for minimal degree rational GBCs.
     Note however that the known proofs of uniqueness for the adjoint requires advanced tools, such as canonical forms.
     We emphasize that our goal is to reconstruct Warren's elementary proof of uniqueness of GBCs without significant external input.
\end{remark}

\subsection{Restriction to facets}
\label{sec:restriction_to_facets}

Later, we will need to perform induction on the dim\-ension of the polytope by stepping from $P$ to its faces.
We therefore need to~under\-stand how rational GBCs (in particular, their degrees) behave under restriction to faces.

If $\bs\alpha=\{p_v/q\}$ are rational GBCs and $F\subset P$ is~a~facet, then the \Def{restriction $\bs\alpha|_F$} yields rational GBCs on $F$:
%
%
%
$$\alpha_v|_F=\frac{p_v|_F}{q|_F}\quad\text{for all $v\in F$}.$$
%
%
We emphasize that this representation is valid only because we have $q|_F\not\equiv 0$ due~to \cref{res:min_not_vanish_of_facet} \itm2.
%
%
This representation might however not be reduced, \ie\ numerator and denominator may still share factors.



\notodo{\msays{Maybe change $r_F$ to $c_F$?}}

\begin{lemma}
    \label{res:r_F}
    Given a facet $F\subset P$ and a vertex $v\in F$ with $\deg(p_v)=m-d$, the following polynomial divides both $p_v|_F$ and $q|_F$:
    $$\text{\Def{$r_F$}}:=\prod \big\{\ell_G|_F \,\big\vert\, \text{$G\subset P$ is a facet neither equal nor adjacent to $F$}\big\}.$$
\end{lemma}
Observe that the shared factor $r_F$ depends neither on the precise GBCs nor the vertex $v\in P$.
This will be important in \cref{res:GBC_deg_m_d_is_WPC}.
\begin{proof}[Proof of \cref{res:r_F}]
    %
    %
    We show that each restricted linear factor $\ell_G|_F$ divides $q|_F$ and all $p_v|_F$. 
    If $G$ is parallel to $F$, then $\ell_G|_F$ is constant, and the claim is trivial.
    If $G$ is not parallel, then since it is not adjacent to $F$, it defines a $(d-2)$-dimensional residual flat $E_G:=\aff(G)\cap\aff(F)$.
    By \cref{res:q_vanishes_residual} and \cref{res:residual_arrangement} $q$ and $p_v$ vanish on $E_G$.
    Since $E_G$ is the zero locus of $\ell_G|_F$ on $\aff(F)$, the claim follows. 
\end{proof}

We can then cancel the common factor $r_F$ provided by \cref{res:r_F} and define~the polynomials $\text{\Def{$p_v^F$}}:=p_v|_F/r_F$ and $\text{\Def{$q^F$}} := q|_F / r_F$.
We call the representation
%
%
$$\alpha_v|_F = \frac{p_v^F}{q^F}\quad\text{for all $v\in F$}$$
the \Def{reduced restriction} of $\{p_v/q\}$ on $F$.
From this definition alone it is not clear~(and in general not true) that we removed all shared factors.
It is however sufficient for preparing the facet induction argument.

\begin{lemma}
    \label{res:degree_on_facet}
    Given a facet $F\subset P$ and a vertex $v\in F$ with $\deg(p_v)=m-d$, then $\deg(p_v^F)= m_F-(d-1)$, where $m_F$ denotes the number of facets of $F$.%
    %
    %
    \notodo{\msays{In this lemma in the proof we put extra effort into establishing that the degree is exactly $m_F-(d-1)$ rather than just upper bounded by this. This cannot be inferred from the known lower bound because for this we would need to cancel at all vertices, but we can only cancel at one vertex because we assume minimal degree only there.}}
    %
    %
\end{lemma}
\begin{proof}
    By definition we have $\deg(p_v^F) = \deg(p_v|_F) - \deg(r_F)$. We are going to~compute the summands on the right individually.

    From \eqref{eq:factorization} we get the factorization
    $$\mathllap{(*)\quad} p_v|_F = \hat p_v|_F \cdot\prod_{G \not\ni v} \ell_G|_F.$$
    From \cref{res:order} \itm1 and \itm2, as well as \cref{res:min_not_vanish_of_facet} we get
    $$o_F(\hat p_v) 
    \overset{\smash{\mathclap{\text{\itm2}}}}= 
    o_F(p_v)
    \overset{\smash{\mathclap{\text{\itm1}}}}=
    o_F(q)
    \overset{\smash{\mathclap{\text{\ref{res:min_not_vanish_of_facet}}}}}=
    0.$$
    Hence, $\hat p_v|_F\not\equiv 0$.
    By \cref{res:hatpv_homogeneous}, $\hat p_v$ is homogeneous at $v\in F$.
    The restriction~of~a homogeneous polynomial to a containing affine subspace is either zero or of the same degree.
    Hence we get $(**)\, \deg(\hat p_v|_F)=\deg(\hat p_v)$. 

    \todo{\msays{I believe somewhere here is the distinction between parallel and not parallel that Warren skipped over. Is this correct? We don't want to make false accusations.}}%
    Assume for now that there is no facet $G$ parallel to $F$.\footnote{This distinction was missing from Warren's proof in \cite{warren2003uniqueness}.}
    Then $(**)\, \deg(\ell_G|_F)=1 = \deg(\ell_G)$ for all facets $G\ne F$. We get
    \begin{align*}
        \deg(p_v|_F) 
        &\,\overset{\mathclap{(*)}}=\, \deg(\hat p_v|_F) + \sum_{F\not\ni v}\deg(\ell_G|_F)
        \\&\,\overset{\smash{\mathclap{(**)}}}=\, \deg(\hat p_v) + \sum_{F\not\ni v}\deg(\ell_G) = \deg(p_v) = m-d.
    \end{align*}
    Likewise
    $\deg(r_F) = m-m_F-1$.
    The claim follows by subtraction.
    
    If, however, $G\ne F$ is a facet that is parallel to $F$, then $\deg(\ell_G|_F)=0$, changing the computed value of both $\deg(p_v|_F)$ and $\deg(r_F)$.
    But both values change~by~exactly 1, and so there is no net change to the value of $\deg(p_v^F)$.
\end{proof}

\begin{corollary}
    \label{res:minimal_to_minimal_of_facet}
    Let $F\subset P$ be a facet.
    If $\{p_v/q\}$ are GBCs on $P$ of minimal~nume\-rator degree $\deg(p_v)=m-d$, then the reduced restrictions $\{p_v^F/q^F\}$ form GBCs~on $F$ of minimal numerator degree $\deg(p_v^F)=m_F-(d-1)$.
\end{corollary}


\begin{remark}\notodo{\msays{Have this here?}}%
    It is this point at which we stop proving per-vertex variants. 
    For~completeness we note that \cref{res:min_deg_o_f_q_eq_ex} and \cref{res:GBC_deg_m_d_is_WPC} below do admit such variants, but their proofs go too far beyond the intended focus of this article. 
\end{remark}

The following result provides an equality analogue of \cref{res:o_q_ge_ex_f}, and a strengthening of \cref{res:min_not_vanish_of_facet} \itm2 to general faces.


\begin{lemma}
    \label{res:min_deg_o_f_q_eq_ex}
    Suppose $\deg(p_v)=m-d$ for all $v$.
    For a face $f\subseteq P$ with $v\in f$~holds
    $$o_f(q) = o_f(p_v) = \ex(f).$$
\end{lemma}
This result is stronger than needed for what follows. Subsequently we only need this statement for when $f$ is a facet or a ridge of $P$.
In both cases $\ex(f)=0$, that~is, \cref{res:min_deg_o_f_q_eq_ex} gives use that $p_v$ and $q$ do not vanish on these faces.
%
%
\begin{proof}[Proof of \cref{res:min_deg_o_f_q_eq_ex}]
    The proof proceeds by induction on $\codim(f)$.
    The statement is clear for $f=P$.
    \notodo{\msays{We could remove the special case for facets because it can be treated by what we see below. If so, we should mention again why $q|_F\not\equiv 0$ etc. But not sure it is maybe cleaner to keep this here.}}If $f$ is a facet then this is precisely \cref{res:min_coprime} \itm2.



    Assume now $\codim(f)\ge 2$.
    We have $o_f(p_v)=o_f(q)\ge\ex(f)$ by \cref{res:order}~\itm2 and \cref{res:o_q_ge_ex_f}.
    Thus, it suffices to show $o_f(q)\le \ex(f)$.
    Fix a facet $F\supset f$.\nls
    Then
    $$o_f(q) \le o_f(q|_F) = o_f(q^F) + o_f(r_F).$$
    We rewrite the summands individually.

    First, we rewrite $o_f(q^F)$.
    \Cref{res:minimal_to_minimal_of_facet} gives us that \{$p_v^F/q^F\}$ form GBCs on $F$ of minima numerator degree $\deg(p_v^F)=m_F-(d-1)$.   
    \notodo{\msays{Careful: true, but only because the polynomials $p_v^F$ and $q^F$ are not zero. This needs to be emphasized somewhere}}
    Also $\codim_F(f)=\codim_P(f) - 1$ (here and in the following we use subscripts to indicate the polytopes with respect to which a quantity is expressed).
    Hence, the induction hypothesis applies and~gives us $o_f(q^F) = \ex_F(f)$
    
    To rewrite $o_f(r_F)$ we use that $r_F$ is the product of restricted linear forms $\ell_G|_F$, where $G\subset P$ is a facet that is neither identical nor adjacent to $F$.
    Such a linear form contributes to $o_f(r_F)$ exactly if $f\subset G$.
    Write $\Fac_P(f)$ (\resp\ $\Fac_F(f)$) for the number of facets of $P$ (\resp\ of $F$) that contain $f$.
    Recall that the facets $G\subset P$ adjacent to $F$ correspond to facets of $F$.
    Then
    \begin{align*}
        o_f(r_F) 
        &=
        \Fac_P(f) - \Fac_F(f) - 1
        \\&=
        (\codim_P(f) + \ex_P(f)) - (\codim_F(f) + \ex_F(f)) - 1
        \\&=
        (\ex_P(f) - \ex_F(f)) + (\codim_P(f) - \codim_F(f) - 1)
        \\&=
        \ex_P(f) - \ex_F(f).
    \end{align*}
    




    \noindent
    This finishes the proof.
\end{proof}

We arrive at the main result of the section (\cf\ \cite[Theorem 6]{warren2003uniqueness}).

\begin{theorem}
    \label{res:GBC_deg_m_d_is_WPC}
    If $\bs\alpha=\{p_v/q\}$ and $\tilbs\alpha=\{\tilde p_v/\tilde q\}$ are rational GBCs with minimal numerator degrees $\deg(p_v)=\deg(\tilde p_v)=m-d$ for all vertices $v\in P$, then $\bs\alpha=\tilbs\alpha$.

    In particular, the Wachspress coordinates are the unique rational GBCs of minimal numerator degree.
\end{theorem}
\begin{proof}
    We proceed by induction on the dimension $d$ of $P$.
    The case $d\le1$ is clear because then $P$ is a simplex on which only barycentric coordinates exist.

    We now assume $d\ge 2$.
    Fix a facet $F\subset P$.
    By \cref{res:minimal_to_minimal_of_facet} the reduced~restrictions $\{p_v^F/q^F\}$ and $\{\tilde p_v^F/\tilde q^F\}$ are GBCs of minimal numerator degree on $F$.
    Hence, the induction hypothesis applies, that is, $\{p_v^F/q^F\}$ and $\{\tilde p_v^F/\tilde q^F\}$ are equivalent and there is a factor $\lambda_F$ so that $p_v^F=\lambda_F \tilde p_v^F$ and $q^F=\lambda_F \tilde q^F$. 
    For all we know, $\lambda_F$ might depend on the facet $F\subset P$. We show that it does not.\footnote{It is this point at which Warren's proof contains a significant gap. Warren does not sufficiently justify why the factors are equivalent across all facets. \Cref{res:r_F} and \cref{res:min_deg_o_f_q_eq_ex} are not proven in \cite{warren2003uniqueness}, but appear necessary in order to complete the argument.}

    From \cref{res:r_F} we get
    \begin{equation}
        \label{eq:lambda_1}
        q|_F 
        = r_F q^F 
        = r_F \lambda_F \tilde q^F 
        = \lambda_F r_F \tilde q^F 
        = \lambda_F \tilde q|_F.
    \end{equation}
    Now choose a ridge $R\subset F$. 
    Restricting both sides of \eqref{eq:lambda_1} to $R$ yields
    $$q|_R = q|_F|_R 
    \overset{\smash{\eqref{eq:lambda_1}}}= \lambda_F \tilde q|_F|_R
    = \lambda_F \tilde q|_R
    \quad\implies\quad 
    \lambda_F = \frac{q|_R}{\tilde q|_R},
    $$
    which is valid since $\tilde q|_R\not\equiv 0$ by \cref{res:min_deg_o_f_q_eq_ex}.
    We see that $\lambda_F$ only depends on the ridge.
    In particular, adjacent facets share the same value for the factor.
    Since the facet-adjacency graph is connected, it follows that $\lambda:=\lambda_F$ is independent of $F$.
    


    At this point we can conclude
    $$r_F p_v|_F=p_v^F = \lambda \tilde p_v^F = \lambda r_F \tilde p_v|_F \quad\implies\quad r_F(p_v - \lambda\tilde p_v)|_F=0.$$
    for all facets $F\subset P$.
    Since $r_F$ does not vanish on $F$, we conclude that each~facet~linear form $\ell_F$ divides $p_v-\lambda\tilde p_v$.
    There are $m$ of these. 
    Since this exceeds the degree~of $p_v$ and $\tilde p_v$, we conclude $p_v=\lambda\tilde p_v$.
    The same argument proves $q=\lambda\tilde q$.
\end{proof}





Induction on \cref{res:minimal_to_minimal_of_facet} then also yields

\begin{corollary}
    \label{res:WPCs_on_facets}
    The restriction of the Wachspress coordinates to some face $f\subset P$ yields the Wachspress coordinates on $f$.
\end{corollary}
%

\subsection{Intermediate degree GBCs}

After proving \cite[Theorem 6]{warren2003uniqueness}, Warren notes
\begin{quote}
    \textit{[\,...\,] }
    the last part of the proof 
    allows for the possibility of rational coordinates functions of degree $n$ \textit{[\,$m$ in our notation\,]} (or greater) that are barycentric. 
    Constructing such a wider class of such coordinates functions is an interesting problem for future research.
\end{quote}
suggesting that numerator degrees strictly between $m-d$ and $m$ are excluded by the same proof.
This is not correct. 
The following example demonstrates that such \Def{intermediate degree GBCs} (\ie\ numerator degree $m-d+k$ with $0<k<d$) do indeed exist and are abundant:

\begin{example}[Intermediate degree]
    Consider the cube $[0,1]^3 = [0,1]^2 \times [0,1]$. 
    The parallelogram has GBCs of degree $4$ (\cf\ \cref{ex:square_perturbation}) and the line has GBCs of degree $1$. 
    The product GBCs on the cube therefore have numerator degree $5 < 6 = m-d-1$. 
    To generalize this we can consider the parallelogram together with a simplex of any dimension. Then we can construct GBCs of degree $5$ for any dimension. 
\end{example}

We do agree that the construction of a wider class of rational GBCs is interesting.
For example, initial investigations (that are beyond the scope of this article) suggest that every polytope that is not a simplex can be equipped with rational GBCs of numerator degree $\deg(p_v)=m-d+k$ for all $k\ge 0$.
It appears interesting and challenging to study the space of rational GBCs on $P$ for a fixed $k$.
For example, the question for $k=1$ seems related to Kohn and Ranestad's \cite[Theorem 7]{kohn2020projective}, which in turn suggests that there is an  $n$-dimensional space of such rational GBCs (here, $n$ is the number of vertices of $P$).

\section{Low denominator degrees}
\label{sec:low_denominator_degrees}

From this section on we transition to the study of denominators of rational~GBCs. 
The focus is on polytopes that admit rational GBCs of especially low denominator degree.
For a polytope $P\subset\RR^d$ let \Def{$\deg(P)$} be the minimally possible denominator degree for rational GBCs on $P$. We call it the \Def{degree of $P$}.
Note that $\deg(P)=0$ if and only if $P$ admits polynomial GBCs, which we address in \cref{sec:polynomial_GBCs}.

Recall that the denominator of the Wachspress coordinates is called the \emph{adjoint polynomial} $\adj_P$ of $P$ and we have $\deg(\adj_P)\le m-d-1$ (\cf\ \cref{sec:degrees}), where $m$ is the number of facets of $P$.
Hence
$$\deg(P) \le \deg(\adj_P)\le m-d-1.$$
A core message of \cite{baumbach2025canonicalformscissorscongruence} is that, in general, $\deg(\adj_P)$ can be much smaller than~$m-d-1$.
The main factor seems to be parallelities in the face structure of $P$, though the precise geometric properties are not well understood.
In particular, and in contrast to numerators, the degree of the Wachspress~denominator does not depend only on the combinatorics of $P$.
It is known that $\deg(\adj_P)\ge m-2d$~\cite[Theorem 3.5 \itm3]{baumbach2025canonicalformscissorscongruence}.

\begin{conjecture}
    For all polytopes $P\subset\RR^d$ holds $\deg(P)\ge m-2d$.
\end{conjecture}

In general, understanding the degree of polytopes seems to be hard. 
We are~aware of the following negative results: 
\begin{enumerate}
    \item There are polytopes where $\det(P)<\det(\adj_P)$. That is, the Wachspress~coordinates are not necessarily the GBCs of minimal denominator degree (\cf\ \cref{ex:adjoint_not_minimal}).
    \item The GBCs of minimal denominator degree $\det(P)$ might not be unique (\cf\ \cref{ex:minimal_not_unique}).
\end{enumerate}
\todo{\msays{Add appendix examples}}

We collect some results about the polytope degree.

\begin{corollary}
    \label{cor:deg_of_faces_does_not_increase}
    If $f\subseteq P$ is a face, then $\deg(f)\le \deg(P)$.
\end{corollary}
\begin{proof}
    It suffices to prove this for when $f$ is a facet, other faces follow by induction.
    If now $\{p_v/q\}$ are rational GBCs on $P$ with $\deg(P)=\deg(q)$, they are in particular maximally canceled. 
    By \cref{res:not_vanish_on_facet} $q$ does not vanish on a facet.
    We can therefore restrict $\{p_v/q\}$ to $f$ and obtain rational GBCs of a denominator degree that cannot be larger than $\deg(q)$.
\end{proof}

Low degree implies combinatorial restrictions:

\begin{lemma}
\label{lem:Shells}
    Let $P_v$ be the vertex figure of $P$ at vertex $v\in P$.
    Every vertex $w\in P_v$ lies in all but at most $\deg(P)+1$ facets of $P_v$.
\end{lemma}




\begin{proof}%
\todo{\msays{Polish this proof.}}
Let $w$ be a vertex of $P_v$. It corresponds to an edge $e=v\bar w$ of $P$, where $\bar w$ denotes the other endpoint of $e$. Let $N$ be the number of facets of $P_v$ that do not contain $w$. We show that $N\leq \deg(P)+1$.

Choose rational GBCs $\{p_z/q\}$ on $P$ such that $\deg(q)=\deg(P)$, and set $r:=o_e(q)$. The point of introducing $r$ is that $q$ may vanish identically on $e$, so we cannot simply restrict the rational functions $p_z/q$ to $e$.

Parametrize $\aff(e)$ by $x_s:=(1-s)v+s\bar w$, and choose a generic vector $u\in\mathbb R^d$ transverse to $\aff(e)$. Since $o_e(q)=r$, we may write
\[
    q(x_s+tu)=t^r\bar q(s)+O(t^{r+1}),
    \qquad
    p_{\bar w}(x_s+tu)=t^r\bar p(s)+O(t^{r+1}),
\]
where $\bar q\in\mathbb R[s]$ is nonzero and $\deg(\bar q)\leq\deg(q)$. The second expansion is possible because \cref{res:order} \itm2 gives $o_e(p_{\bar w})=o_e(q)=r$.

We first relate $\bar p$ and $\bar q$. Choose an affine linear form $L$ with $L(v)=0$ and $L(\bar w)=1$. Affine precision says that $L(x)q(x)=\sum_z L(z)p_z(x)$. Substituting $x=x_s+tu$, we have $L(x_s)=s$, and hence $L(x_s+tu)=s+O(t)$. Thus the coefficient of $t^r$ on the left-hand side is $s\bar q(s)$.

On the right-hand side, the term corresponding to $v$ vanishes because $L(v)=0$, while the term corresponding to $\bar w$ contributes $\bar p(s)$ because $L(\bar w)=1$. If $z\notin e$, then \cref{res:order} \itm2 gives $o_e(p_z)>r$, so these terms do not contribute to the coefficient of $t^r$. Comparing coefficients therefore yields
\[
    \bar p(s)=s\bar q(s).
\]

We now use the factorization (2.3):
\[
    p_{\bar w}
    =
    \hat p_{\bar w}
    \prod_{F\not\ni\bar w}\ell_F.
\]
For every facet $F$ occurring in this product, $\bar w\notin F$, and hence $e\not\subseteq F$. Thus $\ell_F$ does not vanish identically on $e$, so $o_e(\ell_F)=0$. Since $o_e(p_{\bar w})=r$, it follows that $o_e(\hat p_{\bar w})=r$.

We may therefore write $\hat p_{\bar w}(x_s+tu)=t^r h(s)+O(t^{r+1})$ for some $h\in\mathbb R[s]$. Since each $\ell_F$ is affine linear, we also have $\ell_F(x_s+tu)=\ell_F(x_s)+O(t)$. Taking the coefficient of $t^r$ in the factorization above gives
\[
    \bar p(s)
    =
    h(s)\prod_{F\not\ni\bar w}\ell_F(x_s).
\]
From this point on, divisibility is understood in the univariate polynomial ring $\mathbb R[s]$.

Facets of $P_v$ correspond to facets $F$ of $P$ containing $v$, and such a facet of $P_v$ contains $w$ precisely when the corresponding facet $F$ contains $e$. Hence the $N$ facets of $P_v$ not containing $w$ correspond exactly to the facets $F$ satisfying $v\in F$ and $\bar w\notin F$.

For each such facet, $\ell_F(v)=0$ and $\ell_F(\bar w)\neq0$. Since $x_s=(1-s)v+s\bar w$, we have $\ell_F(x_s)=s\ell_F(\bar w)$. Thus each of these $N$ factors contributes a factor $s$ to the product above, and hence $s^N\mid\bar p(s)$. Together with $\bar p(s)=s\bar q(s)$, this gives $s^{N-1}\mid\bar q(s)$ whenever $N\geq1$. Since $\bar q\neq0$,
\[
    N-1\leq\deg(\bar q)\leq\deg(q)=\deg(P).
\]
Thus $N\leq\deg(P)+1$, as claimed.
\end{proof}

\subsection{Polynomial GBCs}
\label{sec:polynomial_GBCs}

In this section we prove the following main result

\begin{theoremX}{\ref{res:polynomial_GBCs}}
    A polytope $P\subset \RR^d$ admits polynomial GBCs if and only if~it~is~af\-finely equivalent to a product of simplices.
\end{theoremX}

As a corollary we obtain

\begin{corollary}
    If $\deg(P)=0$ then $\deg(\adj_P)=0$.
\end{corollary}

As alluded to in \cref{sec:polynomial_GBCs}, polynomial GBCs on cartesian products of simplices are constructed as product coordinates. Polynomiality is moreover preserved under affine transformation.
Thus, it suffices to focus on the ``only if'' part of \cref{res:polynomial_GBCs}.

Throughout this section, let $P\subset\RR^d$ be a polytope with polynomial GBCs.
By the following theorem it suffices to  show that $P$ is simple and affinely equivalent to a 0/1-polytope:

\begin{theorem}[{\cite[Theorem 1]{kaibel2000simple}}]
\label{res:01_polytopes}
    A $0/1$-polytope is simple if and only if~it~is~a~pro\-duct of $0/1$-simplices.
\end{theorem}

That $P$ is simple is a direct consequence of  \cref{lem:Shells}: for denominator degree zero, vertex figures must be simplices:

\begin{corollary}
\label{res:Deg0Simple}
    Polytopes with polynomial GBCs are simple.
\end{corollary}
%

To show that $P$ is affinely equivalent to a 0/1-polytope requires some preparation.
%
%
For this, recall that the \Def{slack matrix} $S_P$ of \(P\) is the matrix whose rows are indexed by the vertices, columns are indexed by the facets, and entries are given by
\[
    (S_P)_{v,F} = \ell_F(v). 
\]
The slack matrix depends on our choice of facet linear forms $\ell_F$, which are determined by $P$ only up to a positive factor.
Hence, we can rescale the columns of $S_P$ freely and still obtain a slack matrix of $P$.
Given the following theorem, it therefore suffices to show that the columns can be rescaled to 0/1-vectors:

\begin{theorem}[{\cite[Theorem 5]{aprile2021recognizing}}]
\label{res:Slack_matrix_affinely_equivalent_to_P}
    \(P\) is affinely equivalent to the convex hull of the rows of its slack matrix $S_P$.
\end{theorem}



With this we have all ingredients for a proof.

\begin{proof}[Proof of \cref{res:polynomial_GBCs}]
Suppose that $P$ has polynomial GBCs $\{p_v\in\RR[x]\}$.

Fix a facet $F\subset P$.
We show that across all vertices $v$, the linear form $\ell_F(v)$ only ever attains two values, depending on whether $v\in F$ or $v\not\in F$. 
Since $\ell_F(v)=0$~for all $v\in F$, the corresponding column of the slack matrix $S_P$ can therefore be rescaled to a 0/1-vector.
The proof then concludes as explained above, using \cref{res:Slack_matrix_affinely_equivalent_to_P}, \cref{res:Deg0Simple} and \cref{res:01_polytopes}.

It remains to show that $\ell_F(v)$ takes on the same value for all $v\not\in F$.
Since~the part of the edge graph of $P$ induced on the vertices outside of $F$ is still connected, it suffices to show that adjacent vertices attain the same value under $\ell_F$.

Suppose then that $v,w\not\in F$ are adjacent via an edge $e\subset P$
and $\ell_F(v)\not=\ell_F(w)$.
Then the restriction $\ell_F|_e$ is not constant.
Moreover, there exists a facet $F'\subset P$~that contains $w$ but not $v$.
Then $\ell_{F'}(w)=0\not=\ell_{F'}(v)$, and the restriction $\ell_{F'}|_e$ is also not constant. 
Restricting the representation \eqref{eq:factorization} to $e$ yields
$$p_v|_e = \tilde p_v|_e \cdot\prod_{G\not\ni v} \ell_G|_e$$
and we see that both $\ell_F|_e$ and $\ell_{F'}|_e$ are independent non-constant factors of $p_v|_e$. 
The latter is therefor at least quadratic.
However, the $p_v|_e$ form GBCs on the~edge, and are therefore necessarily linear.
This is a contradiction, and thus \mbox{$\ell_F(v)=\ell_F(w)$}.

\qedhere

\end{proof}

%% file: sec/appendix.tex
\appendix

\section{Examples}
\label{sec:appendix_examples}
\label{sec:appendix_counterexamples}

We collect here examples to support the main text.

\begin{example}[perturbations]
    \label{ex:perturbations}

    Given rational GBCs $\{\alpha_v\}$, one can construct new rational GBCs $\beta_v:= \alpha_v + r_v$ using a suitable family of rational maps $\{r_v \: P \to\RR\}$.
    Preserving~normalization and precision is equivalent to
    \begin{equation}
        \label{eq:perturbation_conditions}
        \sum_v r_v = 0 \quad\text{and}\quad \sum_v r_v v= 0.
    \end{equation}
    Positivity is harder to guarantee.
    We need at least $\ell_F\mid r_v$ whenever $v\not\in F$ since~we need to preserve the zero loci derived from positivity in \cref{sec:positivity}.
    This is automatic if we instead perturb like
    $$\beta_v := \alpha_v + r_v \prod_F \ell_F.$$

    Consider the square from \cref{ex:square}.
    We can choose $r_{00}=r_{11}=1$ and~$r_{01}=r_{10}=-1$. 
    This yields GBCs of numerator degree 4:
    \begin{align*}
        \beta_{v_1} & = \ell_2 \ell_3 + \ell_1
        \ell_2
        \ell_3
        \ell_4
        , \\
        \beta_{v_2} & = \ell_3 \ell_4 - \ell_1
        \ell_2
        \ell_3
        \ell_4
        , \\
        \beta_{v_3} & = \ell_4 \ell_1 + \ell_1
        \ell_2
        \ell_3
        \ell_4
        , \\
        \beta_{v_4} & = \ell_1 \ell_2 - \ell_1
        \ell_2
        \ell_3
        \ell_4
        .
    \end{align*}
\end{example}

\begin{example}[Adjoint degree is not minimal]
    \label{ex:adjoint_not_minimal}
    We construct a polytope $P$ with $\deg(P)<\deg(\adj_P)$, in dimension $d=3$ with fully explicit GBCs.
    The idea is to start with a polytope with a small adjoint, such as a zonotope, and to move facets parallel to themselves.
    This keeps the normals and parallelisms, but the adjoint jumps back to its generic degree $m-d-1$.
    Other GBCs can still profit from cancellations that the adjoint does not see.

    Let $Z\subset\RR^3$ be the truncated octahedron
    \[
        Z=\{|x|\le6,\ |y|\le6,\ |z|\le6,\ |\pm x\pm y\pm z|\le9\}.
    \]
    It is a zonotope with $m=14$ facets, so $\deg(\adj_Z)=m-2d=8$.
    Let $P$ be obtained from $Z$ by moving three facets parallel to themselves: replace $y\le6$ by $y\le8$, and $x+y+z\le9$ and $x+y-z\le9$ by $x+y+z\le10$ and $x+y-z\le10$ (see \cref{fig:adjoint_not_minimal_3d}).
    Then $P$ is simple with $24$ vertices and the same normals, and $\deg(\adj_P)=10=m-d-1$.
    There are rational GBCs $\{p_v/q\}$ on $P$ with $\deg(q)=9$ and $p_v=\hat p_v\prod_{F\not\ni v}\ell_F$, where each $\hat p_v$ is a quadratic polynomial that is positive on $P$.
    Hence again $\deg(P)\le9<10=\deg(\adj_P)$.
    The exact GBCs are complicatd and are provided on Zenodo. 

\begin{figure}[t]
    \centering
    \includegraphics[width=0.9\textwidth]{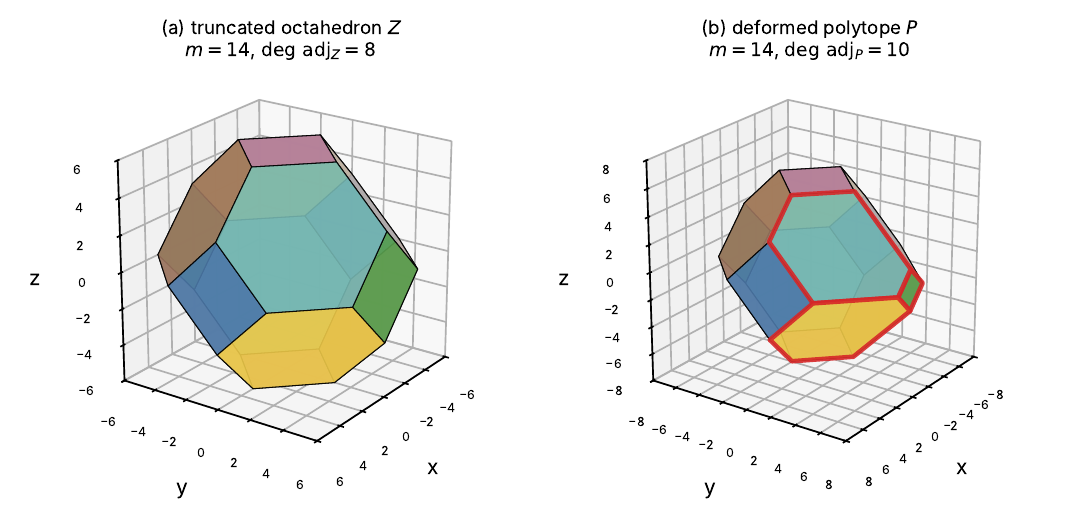}
    \caption{The polytopes of \cref{ex:adjoint_not_minimal} in dimension three. Left: the truncated octahedron $Z$. Right: the polytope $P$ obtained by moving three facets (red edges) parallel to themselves. Both have $14$ facets and the same normals, but $\deg(\adj_Z)=8$ and $\deg(\adj_P)=10$, while $P$ admits rational GBCs with denominator of degree $9$.}
    \label{fig:adjoint_not_minimal_3d}
\end{figure}
\end{example}

\begin{example}[Minimal denominator degree is not unique]
    \label{ex:minimal_not_unique}
    Consider the trapezoid $P=\conv\{v_1,v_2,v_3,v_4\}\subset\RR^2$ with vertices $v_1=(0,0)$, $v_2=(2,0)$, $v_3=(1,1)$, $v_4=(0,1)$.
    It has $m=4$ facets, namely $\{y=0\}$, $\{x+y=2\}$, $\{y=1\}$, $\{x=0\}$, and dimension $n=2$, so $m-n=2$.
    We construct rational GBCs on $P$ whose numerators have the intermediate degree $m-n+1=3$ and the denominator has degree one, but is not equal to the adjoint.
    
    Let
    \begin{align*}
        \p_{v_1}&=(2-x-y)(1-y)(2+x), & \p_{v_2}&=x(1-y)(2+x+y),\\
        \p_{v_3}&=2xy(x+y), & \p_{v_4}&=2y(2-x-y)(1+x),
    \end{align*}
    and $\q=2(2+x-y)$.
    We claim that $\{\p_{v_i}/\q\}$ are rational GBCs on $P$.
    Indeed, one checks by expanding that
    \[
        \sum_{i=1}^4\p_{v_i}=\q,\qquad
        2\p_{v_2}+\p_{v_3}=x\,\q,\qquad
        \p_{v_3}+\p_{v_4}=y\,\q.
    \]
    The first identity is normalization. The other two are linear precision, since $\sum_i\p_{v_i}\,v_i=\big(2\p_{v_2}+\p_{v_3},\ \p_{v_3}+\p_{v_4}\big)$.
    On $P$ we have $0\le x$, $0\le y\le1$ and $x+y\le2$. Hence every factor of every $\p_{v_i}$ is non-negative, and $\q=2(2+x-y)>0$.
    So $\p_{v_i}/\q$ is well defined and non-negative on $P$.

    \emph{Degrees.}
    We have $\deg(\p_{v_i})=3=m-n+1$ for all $i$ and $\deg(\q)=1$.
    The representation is maximally canceled: $\q$ is an irreducible linear form that does not divide $\p_{v_1}=(2-x-y)(1-y)(2+x)$, so $\q$ and the $\p_{v_i}$ have no common non-constant factor.
    Hence the numerator degree $m-n+1$ is not caused by a cancelable common factor of numerators and denominator.

    \emph{Residual flats and the adjoint.}
    The Wachspress denominator of $P$ is $\adj_P=2-y$, which also has degree~$1$.
    However, $\q$ is not a multiple of $\adj_P$.
    The pair of facets $\{x=0\}$ and $\{x+y=2\}$ is not adjacent, so $E=\{(0,2)\}$ is a residual flat, and $\q(0,2)=0$ as required by \cref{res:q_vanishes_residual}; see \cref{fig:residual_arrangement_m_n_plus_1}.
    But $\q$ does not vanish at the point at infinity of the two parallel facets $\{y=0\}$ and $\{y=1\}$:
    in homogeneous coordinates $[x:y:z]$ this point is $[1:0:0]$, and the homogenization $2(2z+x-y)$ of $\q$ takes the value $2\ne0$ there.
    So $\q$ vanishes on the affine, but not on the projective residual arrangement of $P$.
    The adjoint, in contrast, vanishes on the projective residual arrangement.

\begin{figure}
    \centering
    \begin{tikzpicture}[scale=1.5]
        \draw[dotted] (-2.5,0) -- (2.5,0);
        \draw[dotted] (0,-0.5) -- (0,2.5);
        \draw[dotted] (-2.5,1) -- (2.5,1);
        \draw[dotted] (-0.5,2.5) -- (2.5,-0.5);
        \draw[thick,orange!40] (-2.5,-0.5) -- (0.5,2.5);
        \node[left] at (-0.6,1.5) {$\q=0$};
        \fill[blue!45] (0,0) -- (2,0) -- (1,1) -- (0,1) -- cycle;
        \draw[thick] (0,0) -- (2,0) -- (1,1) -- (0,1) -- cycle;
        \draw[->,thick] (-0.8,0.15) -- (-1.2,0.15);
        \node[above] at (-1.0,0.15) {$\infty$};
        \fill[red!50] (0,2) circle (2pt);
    \end{tikzpicture}
    \caption{The trapezoid $P$ of 
    \cref{ex:minimal_not_unique}%
    with its facet lines (dotted), the affine residual flat $(0,2)$ (red point) and the zero set of $\q$ (orange). The parallel facets $\{y=0\}$ and $\{y=1\}$ meet only at infinity.}
    \label{fig:residual_arrangement_m_n_plus_1}
\end{figure}
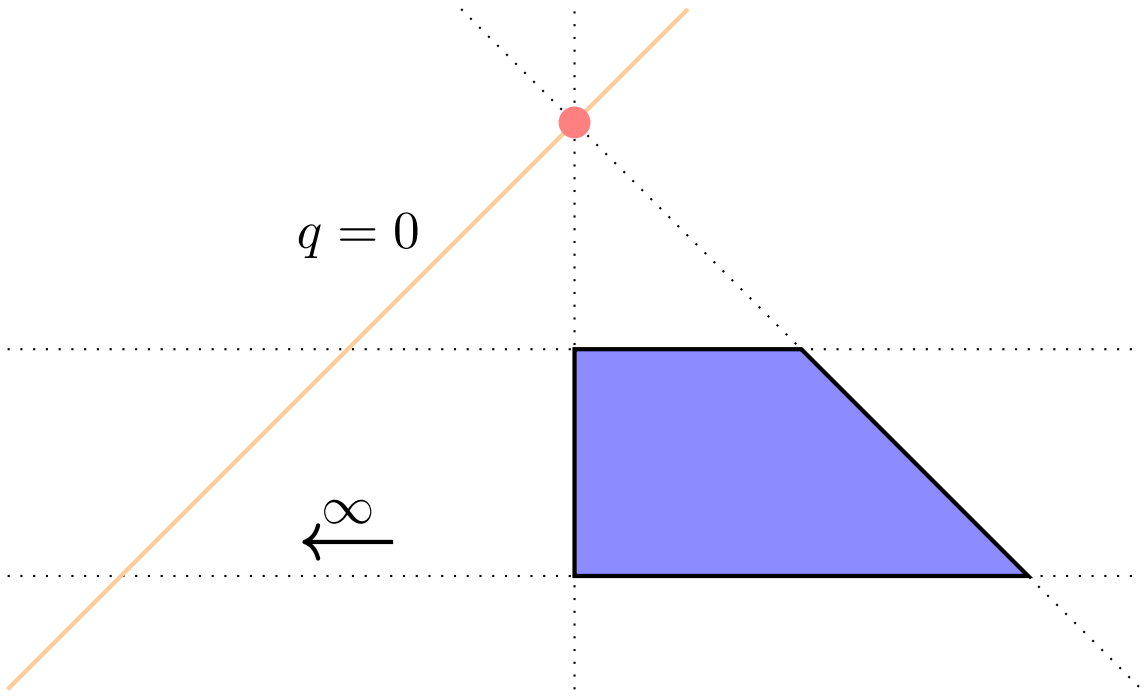
\end{example}

\begin{example}\label{exp:res_arrangement_non_vanishing}
    \cref{res:q_vanishes_residual} states that the denominator $q$ vanishes on residual flats, and \cref{res:residual_arrangement} gives the analogue for a numerator $p_v$ with $v\in E$ under the assumption $\deg(p_v)=m-d$.
    The following example shows that this assumption cannot be dropped.
    There are rational GBCs, a residual flat $E$, and a vertex $v\in E$ with $p_v|_E\not\equiv0$.

\emph{The polytope.}
    Let $Q\subset\RR^3$ be the square pyramid with apex $(0,0,1)$,
    \[
        Q=\{(x,y,z)\in\RR^3:\ 0\le z\le 1,\ |x|\le 1-z,\ |y|\le 1-z\},
    \]
    and let $P=Q\times[0,1]^2\subset\RR^5$ with coordinates $(x,y,z,s,t)$.
    It has $m=5+4=9$ facets, dimension $d=5$, and $20$ vertices, namely $(\sigma,\tau,0,i,j)$ and $(0,0,1,i,j)$ with $\sigma,\tau\in\{\pm1\}$ and $i,j\in\{0,1\}$.
    Thus $m-d=4$.

    \emph{The residual flat.}
    The facets $F=P\cap\{x=1-z\}$ and $\bar F=P\cap\{x=-(1-z)\}$ of $P$ intersect in
    \[
        F\cap\bar F=\{(0,0,1)\}\times[0,1]^2,\qquad
        \operatorname{aff}(F\cap\bar F)=\{x=0,\ y=0,\ z=1\},
    \]
    whereas
    \[
        E:=\operatorname{aff}(F)\cap\operatorname{aff}(\bar F)=\{x=0,\ z=1\}\supsetneq\operatorname{aff}(F\cap\bar F).
    \]
    So $E$ is an intersection of facet-spanning flats that is not spanned by the vertices it contains, \ie\ a residual flat.
    The vertices in $E$ are the four apex vertices $(0,0,1,i,j)$.

    \smallskip
    \noindent\emph{The GBCs.}
    On $Q$ we use the GBCs
    \[
        \alpha_{\mathrm{apex}}=z,\qquad
        \alpha_{(\sigma,\tau)}=\frac{(1-z+\sigma x)(1-z+\tau y)}{4(1-z)},\quad\sigma,\tau\in\{\pm1\},
    \]
    where $(\sigma,\tau)$ denotes the base vertex $(\sigma,\tau,0)$.
    (One checks directly that they are non-negative on $Q$, sum to $1$, and reproduce $x,y,z$.)
    On $[0,1]^2$ we use the bilinear coordinates
    \[
        \lambda_{00}=(1-s)(1-t),\quad \lambda_{10}=s(1-t),\quad \lambda_{11}=st,\quad \lambda_{01}=(1-s)t,
    \]
    and set $B:=\lambda_{00}\lambda_{11}=\lambda_{10}\lambda_{01}=s(1-s)t(1-t)$.
    The product coordinates $\alpha\times\lambda$ of \cref{ex:product} are GBCs on $P$.
    We perturb them at the apex vertices:
    \begin{align*}
        \gamma_{(\sigma,\tau,0,i,j)}&=\alpha_{(\sigma,\tau)}\,\lambda_{ij},\\
        \gamma_{(0,0,1,i,j)}&=z\,\lambda_{ij}+(-1)^{i+j}\,\frac{z\,y^2\,B}{1-z}.
    \end{align*}
    These are rational GBCs on $P$ (verified below) with common denominator $q=4(1-z)$ and numerators
    \begin{align*}
        p_{(\sigma,\tau,0,i,j)}&=(1-z+\sigma x)(1-z+\tau y)\,\lambda_{ij},\\
        p_{(0,0,1,i,j)}&=4z(1-z)\,\lambda_{ij}+(-1)^{i+j}\,4z\,y^2B.
    \end{align*}
    The base numerators have degree $4=m-d$, the apex numerators have degree $7\neq m-d$.

    \emph{Restriction to $E$.}
    Consider the vertex $v=(0,0,1,0,0)\in E$.
    On $E$ we have $x=0$ and $z=1$, so the first summand of $p_v$ vanishes and
    \[
        p_v\big|_E=4\,y^2\,s(1-s)\,t(1-t)\not\equiv0.
    \]
    For instance, at the point $(0,1,1,\tfrac12,\tfrac12)\in E$ it equals $\tfrac14$.
    This point lies in $E$ but not in $P$; on $E\cap P=F\cap\bar F$ we have $y=0$, so $p_v$ does vanish there.
    Since $\deg(p_v)=7\ne m-d$, this does not contradict \cref{res:residual_arrangement}.
    \hfill$\triangle$

\emph{Normalization.}
    Since $\sum_{ij}\lambda_{ij}=1$ and $\sum_{ij}(-1)^{i+j}B=B+B-B-B=0$, the apex coordinates sum to $z$.
    The base coordinates sum to $\sum_{\sigma,\tau}\alpha_{(\sigma,\tau)}=1-z$.
    Hence $\sum_v\gamma_v=1$.

\emph{Linear precision.}
    The perturbation $(-1)^{i+j}\frac{zy^2B}{1-z}$ does not change $\sum_v\gamma_v\,v$:
    it sums to $0$ over $(i,j)$, also when weighted with $i$ (the terms for $(1,0)$ and $(1,1)$ cancel) and with $j$ (those for $(0,1)$ and $(1,1)$ cancel), and the apex has $x=y=0$.
    So $\sum_v\gamma_v\,v$ equals the corresponding sum for the product coordinates, which is $(x,y,z,s,t)$.

\emph{Non-negativity.}
    The coordinates $\gamma_{(\sigma,\tau,0,i,j)}$ are products of non-negative factors.
    On the apex, $\gamma_{(0,0,1,0,0)}$ and $\gamma_{(0,0,1,1,1)}$ are sums of non-negative terms.
    For the remaining two,
    \[
        \gamma_{(0,0,1,1,0)}=z\,\lambda_{10}\Big(1-\frac{y^2\,\lambda_{01}}{1-z}\Big)\ge0,
    \]
    because $\frac{y^2}{1-z}\le1-z\le1$ on $Q$ and $\lambda_{01}\le1$; the vertex $(0,0,1,0,1)$ is analogous.

\emph{Continuity at $z=1$.}
    On $P$ we have $\big|\frac{zy^2B}{1-z}\big|\le(1-z)\,B\le\frac{1-z}{16}\to0$ for $z\to1$.
    Hence $\bs\gamma$ extends continuously to the fiber $\{(0,0,1)\}\times[0,1]^2$ over the apex, where it equals the product coordinates.
    
\end{example}